\documentclass[11pt,letterpaper]{amsart}

\usepackage[latin1]{inputenc}
\usepackage{amsmath,amsthm,amssymb,amscd,mathrsfs}
\usepackage{mathtools,color,esint,bm}
\usepackage{booktabs}
\usepackage[shortlabels]{enumitem}
\usepackage[perpage]{footmisc}
\usepackage{subfig,caption}
\usepackage{placeins}

\definecolor{darkgreen}{rgb}{0.0, 0.6, 0.13}

\usepackage{hyperref}

\newtheorem{thm}{Theorem}[section]

 \theoremstyle{definition}
 
 \theoremstyle{remark}

 \numberwithin{equation}{section}

\newcommand{\Eb}{\mathbb E}

\newcommand{\Rb}{\mathbb R}

\newcommand{\Tb}{\mathbb T}

\newcommand{\Zb}{\mathbb Z}

\newcommand{\Dc}{\mathcal D}

\newcommand {\Jc}{\mathcal{J}}
\newcommand{\Kc}{\mathcal K}
\newcommand{\Lc}{\mathcal L}

\newcommand{\Nc}{\mathcal N}

\newcommand{\Qc}{\mathcal Q}

\newcommand{\Tc}{\mathcal T}

\newcommand{\Ds}{\mathscr D}

\newcommand{\lf}{\mathfrak l}

\newcommand{\nf}{\mathfrak n}

\newcommand{\dirac}{\boldsymbol{\delta}}

\begin{document}
\title{Rigorous Justification of Wave Kinetic Theory}

\author{Yu Deng}
\address{\textsc{Department of Mathematics, University of Southern California, Los Angeles, CA, USA}}
\email{\texttt{yudeng@usc.edu}}
\author{Zaher Hani}
\address{\textsc{Department of Mathematics, University of Michigan, Ann Arbor, MI, USA}}
\email{\texttt{zhani@umich.edu}}
\maketitle
\begin{abstract}
The main purpose of this expository note is to give a short account of the recent developments in mathematical wave kinetic theory. After reviewing the physical theory, we explain the importance of the notion of a ``scaling law", which dictates the relation between the asymptotic parameters as the kinetic limit is taken. This sets some natural limitations on the kinetic approximation that were not precisely understood in the literature as far as we know. We then describe our recent and upcoming works that give the first full, mathematically rigorous, derivation of the wave kinetic theory at the natural kinetic timescale. The key new ingredient is a delicate analysis of the diagrammatic expansion that allows to a) uncover highly elaborate cancellations at arbitrary large order of diagrams, and b) overcome difficulties coming from factorial divergences in the expansion and the criticality of the problem. The results mentioned in this note appear in our recent works \cite{DH21, DH21-2} as well as an upcoming one in \cite{DH22}. 
\end{abstract}

\maketitle
\section{Introduction and the Physical Theory}

The kinetic theory of wave systems, henceforth called Wave Kinetic Theory (WKT), was developed early in the last century, shortly following the advent of Boltzmann's kinetic theory for particle systems. The WKT lies at the center of the theory of non-equilibrium statistical physics for nonlinear waves, which is more commonly known as \emph{wave turbulence theory}. This kinetic theory was developed in parallel in various and distinct fields of physics (firstly by Peierls in the context of anharmonic crystals and Hasselmann in the context of water waves) before it became a systematic paradigm starting from the 1960's with immense practical applications. 

Similar to the particle setting, the aim is to understand the effective behavior of systems with a large number of waves interacting (``colliding") under a weak nonlinear effect. The two main asymptotic parameters in this kinetic theory (at least in the homogeneous setting) are the size $L$ of the domain of the wave system, which measures the number of waves, and the strength $\epsilon$ of the nonlinear interactions. These stand as analogs of the particle number $N$ and particle radius $r$ in Boltzmann's kinetic theory. The kinetic theory emerges under a certain limit of $L \to \infty$ and $\epsilon\to 0$ that is often referred to as the \emph{kinetic limit}. The exact sense that these two limits are taken is specified by the \emph{scaling law} that defines the limiting process, as we shall discuss below.

A longstanding open problem is to put this wave kinetic theory on rigorous mathematical grounds. This means rigorously justifying the approximation of statistical averages of the microscopic dynamics with that of a mesoscopic dynamics given by a kinetic equation known as the \emph{wave kinetic equation}. Here, the microscopic dynamics is given by the dispersive PDE modeling the nonlinear wave system. This is Hilbert's sixth problem for nonlinear waves. It should be pointed out that the importance of this problem is not merely mathematical, but as we shall discuss below, has a direct impact on the physical theory and on the proper formulation of the wave kinetic theory, its setup, regimes, expectations, and conclusions. 

\subsection{Physical theory} Let us start by reviewing the physical theory. Weak nonlinear interactions yield a wide scale separation between linear timescales and nonlinear ones. The angles of the Fourier modes $A_k$ at wave number $k$ evolve on the fast linear timescale, whereas the amplitudes $|A_k|$ evolve on much slower ones. The average distribution of the angles on longer timescales, at least in absence of nonlinear interaction, is given by the uniform measure on the unit circle. For this reason, a \emph{random phase} (RP) assumption is often assumed (with various degrees of rigor) in the formulation of the wave kinetic theory. Of course, in a rigorous justification like the one we consider here, one can only adopt such an assumption at initial time, and would have to justify an appropriate version of how it propagates at later times. In many ways, understanding this so-called \emph{propagation of chaos} from a qualitative and quantitative perspective, is the heart of the problem. Given the fast-slow scale separation between the dynamics of the angles and amplitudes of $A_k$, the RP assumption is natural from the point of view of the physical theory, and is indicative of the generic behavior of solutions with deterministic data.

As the phases get averaged out, the central quantity in WKT are the amplitudes $|A_k|$ of the waves oscillating at wave number $k$. The variance of these amplitudes is described, in the above kinetic limit, by the \emph{wave kinetic equation} (WKE); here the averaging is done with respect to the (random phase) distribution of Fourier modes at initial time. The form of the WKE depends on the microscopic dynamics that govern the wave interactions (including the dispersion relation $\omega_k$ and the nonlinearity). To make the discussion concrete, let us fix our emblematic microscopic wave system given by the nonlinear Schr\"odinger (NLS) equation
\begin{equation}\label{NLS}\tag{NLS}
(i\partial_t +\Delta)u=\epsilon |u|^2 u, \qquad t\in \Rb, x\in \Tb^d_L
\end{equation}
where the spatial domain $\Tb^d_L$ is a periodic box of size $L$, and $\epsilon>0$ denotes the strength of the nonlinear interaction. In Fourier space, this equation takes the form
\begin{equation}\label{Feqn}
i\partial_t A_k -\omega(k) A_k=\epsilon L^{-d} \sum_{(k_1, k_2, k_3)\in (\Zb_L^d)^3} \boldsymbol{T_{k_1 k_3}^{k_2 k} }A_{k_1}\overline{A_{k_2}}A_{k_3}.
\end{equation}
Here $k=(k^{1}, \ldots, k^{d}) \in \Zb^d_L:=L^{-1}\Zb^d$, which is a mesh of spacing $L^{-1}$ and tends to continuum in the limit $L\to \infty$. The Fourier coefficient at wave number $k$ is given by $A_k=L^{-d/2} \int u(x) e^{-i k \cdot x}\,\mathrm{d}x$, and $\omega(k)=|k|_\beta^2:=\sum_{j=1}^d\beta_j (k^{j})^2$ is the dispersion relation which depends on the aspect ratios $\beta:=(\beta_1,\cdots,\beta_d)$ of the torus (e.g. $\beta_j=1$ in the case of square torus). Finally, $\boldsymbol{T_{k_1 k_3}^{k_2 k}}=T_{k_1 k_3}^{k_2 k}\delta(k_1+k_3-k_2-k)$ with $T_{k_1 k_3}^{k_2 k}=1$. In fact, equation \eqref{Feqn} represents the general form for a nonlinear wave system with cubic interactions and dispersion relation $\omega(k)$. We shall refer to this general form later in our discussion.

\medskip 

As mentioned earlier, the initial data is taken to satisfy the random phase assumption as follows: 
\begin{equation}\label{data}\tag{DATA}
A_k(t=0)= (A_{\mathrm{in}})_k:=\sqrt {n_{\mathrm{in}}(k)} \cdot \eta_k
\end{equation}
where ${n_{\mathrm{in}}(k)}$ is a non-negative, sufficiently smooth, and decaying function on $\Rb^d$ which represents the initial amplitudes, and $\eta_k$ are i.i.d. complex random variables with mean zero and variance 1 (normalized). In particular, we have that $\Eb |A_k(t=0)|^2=n_{\mathrm{in}}(k)$. We shall impose a few more assumptions on the law of the $\eta_k$ to state our theorems, but they include standard choices of interest like Gaussians or random phase (i.e. $\eta_k$ uniformly distributed on the unit circle $\mathbb S^1$). We emphasize that this independence assumption of distinct $A_k$ is only adopted here at initial time. 

The WKT states that, with $A_k(t)$ satisfying \eqref{Feqn} and \eqref{data}, the evolution of $\Eb |A_k(t)|^2$ is approximated, in the limit of asymptotically large $L$ and small $\epsilon$, by solutions of the wave kinetic equation
\begin{equation}\label{wke}\tag{WKE}
\left\{
\begin{aligned}&\partial_t n(t,k)=\Kc(n(t))(k),\\
&n(0,k)=n_{\mathrm{in}}(k),
\end{aligned}
\right.
\end{equation} where the nonlinearity
\begin{multline}\label{wkenon}\tag{KIN}\Kc(\phi)(k)=\int_{(\Rb^d)^3}|T_{k_1 k_3}^{k_2 k}|^2\phi(k)\phi(k_1)\phi(k_2)\phi(k_3)\bigg(\frac{1}{\phi(k)}-\frac{1}{\phi(k_1)}+\frac{1}{\phi(k_2)}-\frac{1}{\phi(k_3)}\bigg)\\\times \dirac_{\Rb^d}(k_1-k_2+k_3-k) \cdot \dirac_\Rb(\omega(k_1)-\omega(k_2)+\omega(k_3)-\omega(k))\,\mathrm{d}k_1\mathrm{d}k_2\mathrm{d}k_3.
\end{multline}
Here $\dirac_{\Rb^d}$ is the Dirac delta on $\Rb^d$, and in the case of (\ref{NLS}), $|T_{k_1 k_3}^{k_2 k}|^2=1$ and $\omega(k)=|k|_\beta^2$. 

This approximation happens at the \emph{kinetic (or Van-Hove) timescale} $T_{\mathrm{kin}}:=\epsilon^{-2}$ in the sense that
\begin{equation}\label{approx}
\Eb |A_k(T_{\mathrm{kin}}\cdot \tau)|^2 \approx n(\tau, k), \qquad \textrm{ as }L\to \infty \textrm{ and } \epsilon \to 0.
\end{equation}

A formal derivation\footnote{Such formal derivations can be found in different guises in various physical treatments \cite{ZLFBook, Nazarenko}. We present it here in a fashion closer to our mathematical understanding and in a way that allows to point out some qualitative features.} of \eqref{wke} from \eqref{Feqn} proceeds via an expansion of the solution to \eqref{Feqn} into its first Picard iterates (via a multiple timescale expansion). Under the RP assumption, the first order term (in $\epsilon$) in the corresponding expansion of $\Eb |A_k|^2$ will vanish. The second order terms that are $O(\epsilon^2)$ give a closed form in $n_{\mathrm{in}}(k)=\Eb |(A_{\mathrm{in}})_k|^2$ and can be written as a sum of four terms, one of which is given by
\begin{equation}\label{formal1}
\epsilon^2 t^2L^{-2d}\sum_{(k_1, k_2, k_3)}|T_{k_1 k_3}^{k_2 k}|^2 n_{\mathrm{in}}(k_1)n_{\mathrm{in}}(k_2)n_{\mathrm{in}}(k_3) \dirac_{\Rb^d}(k_1-k_2+k_3-k) \times \left(\frac{\sin(\pi \Omega t)}{\pi \Omega t}\right)^2,
\end{equation}
where $\Omega=\omega(k_1)-\omega(k_2)+\omega(k_3)-\omega(k)$. This sum can be formally approximated\footnote{Depending on the adopted scaling law (cf. Section \ref{sec:scalinglaws}), the rigorous proof of this statement can often require deep results in analytic number theory, as we shall see later.} in the limit of large $L$ by the integral
\begin{equation}\label{intlim}
\epsilon^2 t^2\int_{(\Rb^d)^3}|T_{k_1 k_3}^{k_2 k}|^2 n_{\mathrm{in}}(k_1)n_{\mathrm{in}}(k_2)n_{\mathrm{in}}(k_3)  \dirac_{\Rb^d}(k_1-k_2+k_3-k) \times \left(\frac{\sin(\pi \Omega t)}{\pi \Omega t}\right)^2\,\mathrm{d}k_1\mathrm{d}k_2\mathrm{d}k_3.
\end{equation}
Finally, notice that $t \left(\frac{\sin(\pi \Omega t)}{\pi \Omega t}\right)^2 $ is an approximation  identity which converges in the limit of large $t$ to $\dirac(\Omega)$. The limit of (\ref{intlim}) then equals $\epsilon^2 t$ multiplied by the first term in \eqref{wkenon}, but with $\phi$ replaced by $n_{\mathrm{in}}(k)$. The other terms that are $O(\epsilon^2)$ give the other terms in the first iterate of \eqref{wke}. This agreement at the first few iterates between $\Eb |A_k|^2$ and $n(t, k)$ in the limit of large $L$ and small $\epsilon$ (which guarantees the largeness of $t\sim T_{\mathrm{kin}}\sim \epsilon^{-2}$ in the approximation identity above) leads one to conjecture the validity of the approximation in \eqref{approx}.
 
\section{Scaling laws in wave kinetic theory}\label{sec:scalinglaws}
\subsection{General Considerations} What is missing from the above description is the precise manner in which the two limits $L \to \infty$ and $\epsilon\to 0$ are taken. This is the notion of the \emph{scaling law} which specifies precisely the relation between $\epsilon$ and $L$ in the kinetic limit. Such scaling laws are not uncommon in kinetic theories; in fact, the derivation of the Boltzmann kinetic equation for particle systems \cite{Lanford, GSRT} is done in the so-called Boltzmann-Grad scaling law \cite{Grad} for which $N r^{d-1}=O(1)$ as the particle number $N\to \infty$ and particle radius $r\to 0$. It is safe to say that a systematic understanding of the role and nuance of scaling laws in wave kinetic theory had not been fully developed prior to the recent rigorous mathematical investigations, at least to the best of our knowledge. In fact, the ramifications of these investigations and their implications on the physical theory raise several intriguing questions.
\medskip

The most general form of a scaling law between $\epsilon$ and $L$ is to set $\epsilon=L^{-\gamma}$ where $0\leq \gamma \leq \infty$, where by $\gamma=0$ we mean that the limit $L\to \infty$ is taken first followed by the limit $\epsilon \to 0$, and vice versa for the case $\gamma=\infty$. As we shall promptly see, not all scaling laws, now characterized by the value of $\gamma\in [0, \infty]$, are admissible for the kinetic theory. To see this, one can just look at the formal derivation we have given in the previous section. According to the heuristic derivation presented above, the sum in \eqref{formal1} has size $O(L^{2d}t^{-1})$ which is coming from a highly nontrivial estimate on the number of $\Zb_L^d$ lattice points in the window $W_{t^{-1}}\approx W_{\epsilon^2}$ (since $t\sim T_{\mathrm{kin}}=\epsilon^{-2}$) around the resonant manifold $\{\Omega=0\}$, where
$$
W_\delta=\{(k_1, k_2, k_3): k_1-k_2+k_3=k, |\Omega|\leq \delta, \textrm{ and }|k_j|\lesssim 1\}.
$$
The last restriction on $k_j$ comes from the decay properties of $n_{\mathrm{in}}$. In fact, this estimate, which compares the number of such lattice points in $W_{t^{-1}}$ with $L^{2d}\times\mathrm{Vol}(W_{t^{-1}})$, is not true unconditionally for very large $t$. An easy way to see this limitation is to notice that if all the $\beta_j=1$, then the dispersion relation $\omega(k)$ takes values in $L^{-2}\Zb$, so if $t>L^{2}$, the number of lattice points in $W_{t^{-1}}$ would be constant (equals that of $W_0$) and independent of $t$.

This upper bound on $t$ is what gives the limitation on the admissible scaling laws. In fact, at the the kinetic time $T_{\mathrm{kin}}=\epsilon^{-2}=L^{2\gamma}$, such an upper bound translates into an upper bound $\gamma_{\mathrm{max}}$ on the admissible $\gamma$. The exact value of $\gamma_{\mathrm{max}}$ depends on the dispersion relation $\omega(k)$ and lattice-point equidistribution properties associated to it. This is often a deep question in analytic number theory, and at this point, we only have a good understanding of the answer in the NLS case, where $\omega(k)=|k|_\beta^2$. Using results on equidistribution of rational and irrational quadratic forms \cite{FGH, BGHS1, BGHS2}, we know that for the square torus $\beta=(1, 1, \ldots, 1)$ we have $\gamma_{\mathrm{max}}=1$, whereas for generic irrational\footnote{This means that the $\beta_1, \ldots, \beta_d$ are rationally independent and skip a measure zero set.} $\beta$, this $\gamma_{\mathrm{max}}$ can be as large as $d/2$.

\medskip 

For more general $\omega(k)$ where the value of $\gamma_{\mathrm{max}}$ might be hard to find, one can identify \emph{sufficient} conditions on $\gamma$ that guarantee the needed equidistribution result. In effect, this gives a lower bound on $\gamma_{\mathrm{max}}$. This can often be done without the need to resort to sophisticated number theoretic techniques. Heuristically, it can be stated as the requirement that the minimum separation between values of $\omega(k)$ with nearby $k$ is much less than the width $\delta$ of the window $W_\delta$. This leads to the condition \cite{Nazarenko, GBE22}
$$
\frac{1}{L} \left|\frac{\partial \omega}{\partial k}\right| \ll \frac{1}{T_{\mathrm{kin}}}\sim L^{-2\gamma},
$$
which gives the upper bound $\gamma\le 1/2$. This upper bound is sufficient but is not necessary for most dispersion relations. For NLS, it is not necessary regardless of the aspect ratio of the torus (square or irrational), as mentioned above. That being said, such upper bounds can be the best ones available in cases when no number theoretic techniques are known.

\begin{figure}[h!]
  \includegraphics[scale=0.29]{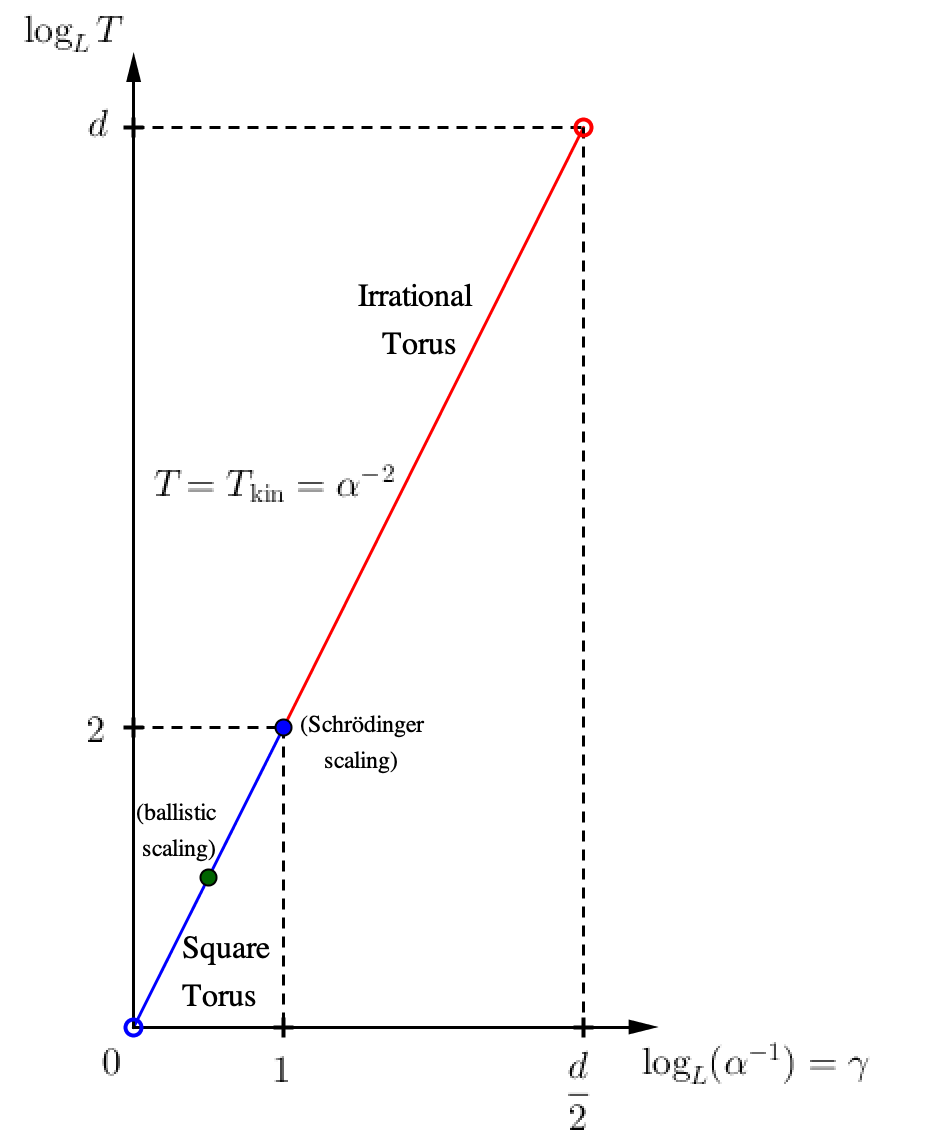}
  \caption{Log-Log plot of scaling laws, represented by $\gamma$, versus time, represented by $\log_L T$. The line represents the kinetic timescale $T=\alpha^{-2}$. The blue segment corresponds to the scaling laws $0<\gamma<1$, which is the admissible range for arbitrary (in particular square) tori and is covered by Theorem \ref{main}. The full segment corresponds to the scaling law $0<\gamma<d/2$, which is admissible for generic irrational tori. The two specific scaling laws ($\gamma=1$ and $\gamma=1/2$, see Section \ref{specscaling}) are also indicated.}
  \label{fig:SL}
\end{figure} 
  \subsection{Exact resonances, quasi resonances, and scaling laws} Another way to frame the discussion on the admissible scaling laws is in terms of the comparison between exact and quasi-resonances. Exact resonances correspond to the interactions for which $\Omega=0$, whereas quasi-resonances are ones that correspond to interactions such that $0<|\Omega|\lesssim t^{-1}$. While exact resonances play an important role in the study of dispersive PDE as deterministic dynamical systems, it is clear from the formal derivation presented above that the nonlinearity in \eqref{wke} comes as the thermodynamic limit of the quasi-resonant interactions. For this kinetic term to be the leading contribution, one needs to have the quasi-resonant interactions dominate the exact resonant ones. The size of the exact resonant interactions depends on the dispersion relation $\omega(k)$, and often also requires some input from number theory. 

In the case of NLS, the contribution of the exact resonances (obtained by replacing $\left(\frac{\sin(\pi \Omega t)}{\pi \Omega t}\right)^2$ with $\mathbf{1}_{\Omega=0}$) to \eqref{formal1} is $\epsilon^2 t^2 L^{-2}$ for the square torus case $\beta=(1, \ldots, 1)$, and $\epsilon^2t^2 L^{-d}$ for generic irrational $\beta$ \cite{BGHS1, BGHS2}. Comparing this to the size $\epsilon^2 t$ of the quasi-resonances (which is given formally by the volume counting argument mentioned above), one obtains that the quasi-resonances dominate the exact ones precisely when $t<L^2$ in the square torus case, and when $t< L^d$ in the rationally independent case. At the kinetic time $t\sim \epsilon^{-2}=L^{2\gamma}$, we again obtain {\bf the requirements that $\gamma < 1$ for the square torus, and $\gamma <  d/2$ for the irrational torus}, which coincide with the $\gamma_{\mathrm{max}}$ values mentioned in the previous subsection. This shows that if one considers scaling laws $\epsilon=L^{-\gamma}$ with $\gamma>\gamma_{\mathrm{max}}$, then the limit dynamics will be given by the exact resonances, and not the kinetic equation. This direction was investigated in \cite{FGH, BGHS1} leading to the so-called Continuous Resonant (CR) equations. 
\subsection{Significance of certain scaling laws}\label{specscaling} For the rest of this note, we will restrict the discussion to (\ref{NLS}). In this setting, a couple of scaling laws seem to be of particular importance. The first one corresponds to the case $\gamma=1$, so that $\epsilon=L^{-1}$ and $T_{\mathrm{kin}}=L^2$. The importance of this scaling law comes from the natural parabolic scaling of the NLS equation, by which the solutions to \eqref{NLS} on the torus of size $L$ and at the timescales $\sim L^2$ can be rescaled to solutions to the cubic NLS equation on the unit torus and at $O(1)$ times. Namely,
$$
\textrm{If }v(t,x)= L^{\frac{1}{2}} u(L^2 t, Lx), \qquad  x\in \Tb^d, \qquad  \textrm{then}\qquad (i\partial_t v+\Delta)v=|v|^2v.
$$
This means that the conclusions of the kinetic theory under this scaling law, can be translated into conclusions about the dynamics of NLS on the unit torus. In three dimensions, this is closely related to the famous Gibbs measure invariance problem for cubic NLS (i.e. invariance of the $\Phi_3^4$ measure under the Schr\"{o}dinger dynamics), which is the only Gibbs measure invariance problem that still remains open in all polynomial settings after the works of \cite{Bourgain94,Bourgain96,Zhidkov,OT,DNY,BDNY}. In addition, energy cascade behavior for NLS can also be observed at the level of (\ref{wke}) \cite{Nazarenko, EV}, and proving such cascade dynamics for the NLS equation on the unit torus is a problem of great interest \cite{Bourgain}.

Another important scaling law for (\ref{NLS}) is $\gamma=\frac12$ for which $T_{\mathrm{kin}}=L$. This scaling law can be called the ballistic scaling law because it equates the kinetic timescale $\epsilon^{-2}$ with the ballistic timescale needed for a wave packet at frequency $O(1)$ to traverse the domain $\Tb^d_L$. This seems analogous, as we shall explain, to the Boltzmann-Grad scaling law adopted in the derivation of Boltzmann's equation, in which the so-called mean-free path is also equated to the transport length scale. It should be pointed though that such wave packet considerations are more relevant in the inhomogeneous setting of the problem, where the initial field is not assumed to be homogeneous in space as in \eqref{data}. An example of such data is when one sets \eqref{NLS} on $\Rb^d$ and adopts data of the following form. Fix a function $W_0=W_0(x,\xi):\mathbb{R}^3\times\mathbb{R}^3\to\mathbb{R}_{\geq 0}$ which decays rapidly in $\xi$ and $x$. Consider the distribution of random initial data $u_{\mathrm{in}}$ of \eqref{NLS}, such that the Wigner transform
\begin{equation}\label{wigner}\mathbb{E}\bigg(\int_{\mathbb{R}^3}e^{iLy\cdot\eta}\,\overline{\widehat{u_{\mathrm{in}}}\big(\xi-\frac{\eta}{2}\big)}\widehat{u_{\mathrm{in}}}\big(\xi+\frac{\eta}{2}\big)\,\mathrm{d}\eta\bigg)\to W_0(y,\xi)\quad (\mathrm{as\ }L\to\infty),
\end{equation} (possibly in a weak sense). This is achieved, for example, by setting the random data as
\begin{equation}\label{data0}
u_{\mathrm{in}}(x)=L^{-\frac{d}{2}}\sum_{k\in(L^{-1}\mathbb{Z})^3}\psi\bigg(\frac{x}{L},k\bigg)\cdot g_k\cdot e^{ik\cdot x};\quad \psi(y,k)=\sqrt{W_0(y,k)},
\end{equation} which can be viewed as an inhomogeneous generalization of that in \eqref{data}. Then, the solution to (\ref{NLS}) has the form
\begin{equation}\label{solution}
u(t,x)=L^{-\frac{d}{2}}\sum_{k\in(L^{-1}\mathbb{Z})^3}A\bigg(t,\frac{x}{L},k\bigg)\cdot e^{ik\cdot x}.
\end{equation} Denoting $N(t,y,k):=\mathbb{E}|A(t,y,k)|^2$, which corresponds to the Wigner transform of $u(t)$, and performing a formal expansion similar to the one in \eqref{formal1}, we find that $N$ satisfies
\begin{equation}\label{wke1}\partial_tN+\frac{1}{L}(k\cdot\nabla_y)N=\epsilon^2 \mathcal{C}(N,N,N)+\mathrm{l.o.t.}, \qquad N(0,y,k)=W_0(y,k).
\end{equation}
This gives the inhomogeneous wave kinetic equation provided one equates the transport timescale $L$ with the kinetic timescale $\epsilon^{-2}$, which is the scaling law $\gamma=\frac12$ with $T_{\mathrm{kin}}=L$. 

It is worth pointing out that, if one would like to see the homogeneous WKE as a limit of the inhomogeneous one, then one has to introduce an additional parameter to the data in \eqref{data0}, namely one measuring the scale of the inhomogeneity. This can be done by rescaling $W_0$, or equivalently by replacing $\psi(\frac{x}{L},k)$ with $\psi(\frac{x}{M},k)$ in \eqref{data0}, where $M$ is the new inhomogeneity parameter. This leads to the flexibility of scaling laws in the homogeneous setting; in fact all the admissible scaling laws described above arise as suitable limits with $L\to\infty$ and $M/L\to \infty$.
\section{Rigorous Justification Theorems} We are now ready to state the rigorous theorems justifying the approximation (\ref{approx}) at the kinetic timescales. These are proven in \cite{DH21, DH21-2, DH22}, which follow several related works and partial results in \cite{EY, ESY, LukSpohn, BGHS2, DH19, CG1, CG2, DK1, DK2, DKMV} (see also \cite{ACG, DGHG, Ma} for more recent related developments, particularly in the stochastically forced setting of \cite{ST}).

\begin{thm}[Deng-Hani 2021]\label{main} Consider (\ref{NLS}) on the periodic box $\Tb^d_L$ with $d\geq 3$.
\begin{itemize}
\item Take $n_{\mathrm{in}}\geq 0$ sufficiently regular and decaying, and the initial data $u_{\mathrm{in}}$ to be well-prepared with Fourier coefficients $(A_{\mathrm{in}})_k$ given in \eqref{data}, and suppose that the law of $\eta_k$ is rotationally symmetric and has exponential tails (e.g. Gaussian). 
\item {\it Scaling laws:} Let $\alpha\sim L^{-\gamma}$ for $\gamma \in (0,1]$, and recall that $T_{\mathrm{kin}}=\alpha^{-2}$. For $\gamma=1$, we assume suitable genericity conditions on the aspect ratios $\beta$ of the box. 

\end{itemize}

THEN, 
\begin{enumerate}
\item \emph{Justification of WKE:} There exists $\delta<1$ \emph{fixed}, and an absolute constant $\nu>0$ such that, for $L$ large enough, it holds that 
$$
\mathbb E \left|A_k(\tau \cdot T_{\mathrm{kin}})\right|^2 = n(\tau, k)+O(L^{-\nu})
$$
uniformly in $(\tau, k) \in [0, \delta]\times \Zb^d_L$. Here $n(\tau,k)$ solves the wave kinetic equation (\ref{wke}) with data $n_{\mathrm{in}}$.

\medskip

\item \emph{Propagation of Chaos:} If $k_1, \ldots, k_r$ are distinct wave numbers, the random variables $A_{k_j}(t)$ $(1\leq j \leq r)$ retain their independence in the kinetic limit $L\to \infty$ in the following sense:
\begin{equation}\label{moments}
\Eb\big(A_{k_1}^{p_1}(t)\ldots A_{k_r}^{p_r}(t)\cdot\overline{A_{k_1}^{q_1}(t)}\ldots \overline{A_{k_r}^{q_r}(t)}\big)=\prod_{j=1}^r \delta_{p_j=q_j}\, \Eb |A_{k_j}(t)|^{2p_j}+O(L^{-\nu})
\end{equation}
 for {$t\in [0, \delta \cdot T_{\mathrm{kin}}]$}. Moreover, the asymptotics of $\Eb |A_{k_j}(t)|^{2p_j}$ can be completely described in terms of the solution $n(\frac{t}{T_{\mathrm{kin}}}, k)$ of (\ref{wke}). 
  
 \medskip

\item \emph{Limiting law:} The law of $A_k(\tau \cdot T_{\mathrm{kin}})$ converges to the law with density function $\rho_k(\tau, v)$ (with $v\in \Rb^2$) which evolves according to the linear PDE
$$
\partial_\tau \rho_k=\frac{\sigma_k(\tau)}{4} \Delta \rho_k -\frac{\gamma_k(\tau)}{2} \nabla \cdot(v \rho_k),
$$
where $\sigma_k(\tau)>0$ and $\gamma_k(\tau)$ are functions explicitly defined in terms of the solution $n(\tau,k)$ to the wave kinetic equation. 

\item \emph{Propagation of Gaussianity:} In particular, if each $\eta_k$ is Gaussian, then $\rho_k(\tau, v)$ is Gaussian with variance $n(\tau, k)$ for any $\tau>0$.

\end{enumerate}

\end{thm}

This theorem provides the rigorous mathematical foundation of wave kinetic theory after almost a century of its inception as a physical theory. It clarifies for the first time, the subtle role played by scaling laws, and provides the needed delicate analysis of the Feynman diagram expansion that underlies this theory. Indeed, as we shall explain below, this analysis is heavily dependent on the scaling laws, as certain non-leading diagrams can be small error terms for certain scaling laws, and can be too large and divergent for others, thus requiring uncovering subtle cancellations amongst them. While the first two conclusions in Theorem \ref{main} were widely conjectured in all the physics literature on wave kinetic theory, it should be pointed out that the third one has been less well-known. In a way, it follows as a consequence of Peierls' original treatise \cite{Peierls}, but it was more recently discovered in the form above in the physics literature in \cite{CLN,CJKN}. As such, the above rigorous results again confirm the physical predictions.

\medskip

Before we move on to discuss the diagrammatic expansion and the elaborate cancellations involved, let us remark that the special moments in \eqref{moments} with $p_j=q_j=1$ converge, at the kinetic timescale $t=\tau \cdot T_{\mathrm{kin}}$, to $n_{r}(\tau, k_1, \ldots, k_r)=\prod_{j=1}^r n(\tau, k_j)$ which is a solution of the wave kinetic hierarchy (WKH) given by
\begin{equation*}
\begin{split}
\partial_t n_r(t, k_1, \cdots, k_r)&=\sum_{j=1}^r\int_{(\Rb^d)^3}\dirac(\ell_1-\ell_2+\ell_3-k_j)\cdot \dirac(|\ell_1|_\beta^2-|\ell_2|_\beta^2+|\ell_3|_\beta^2-|k_j|_\beta^2)\,\mathrm{d}\ell_1\mathrm{d}\ell_2\mathrm{d}\ell_3\\
&\!\!\!\!\!\!\!\!\!\!\!\!\!\!\!\!\!\!\!\!\!\!\!\!\!\!\!\!\!\!\times \bigg[ n_{r+2}(t, k_1, \cdots, k_{j-1}, \ell_1, \ell_2, \ell_3, k_{j+1},\cdots, k_r)+n_{r+2}(t, k_1, \ldots, k_{j-1}, \ell_1, k_j, \ell_3, k_{j+1},\cdots, k_r)\\
& \!\!\!\!\!\!\!\!\!\!\!\!\!\!\!\!\!\!\!\!\!\!\!\!\!\!\!\!\!\!-n_{r+2}(t, k_1, \cdots, k_{j-1}, k_j, \ell_2, \ell_3, k_{j+1},\cdots, k_r)-n_{r+2}(t, k_1, \cdots, k_{j-1}, \ell_1, \ell_2, k_j ,k_{j+1},\cdots, k_r)\bigg].
\end{split}
\end{equation*} 
This hierarchy stands as an analog to Boltzmann and Gross-Pitaevskii hierarchies, and was formally derived in recent works such as Chibarro {et al.} \cite{CDJR1,CDJR2}, Eyink-Shi \cite{EShi} and Newell-Nazarenko-Biven \cite{NNB}, though it also follows from earlier works including the foundational work of Peierls, see \cite{Peierls, BP, Nazarenko}.

The property that factorized initial data of form $(n_r)_{\mathrm{in}}(k_1, \cdots, k_r)=\prod_{j=1}^r n_{{\mathrm{in}}} (k_j)$ lead to factorized solutions of the WKH of form $n_r(\tau, k_1, \ldots, k_r)=\prod_{j=1}^r n(\tau,k_j)$ where $n(\tau,k)$ solves (\ref{wke}) with initial data $n_{\mathrm{in}}$, is called \emph{factorizability}. As such, Theorem \ref{main} gives directly a rigorous derivation of the WKH for factorized solutions. In \cite{DH21-2}, we actually give an essentially full rigorous derivation of the WKH for more general, non-factorized solutions, by applying the quantum de Finetti theorem and representing the general case as superposition of factorized cases. We refer the reader to \cite{DH21-2} for details.

\section{Diagrammatic Expansion}

\begin{figure}[h!]
  \includegraphics[scale=0.41]{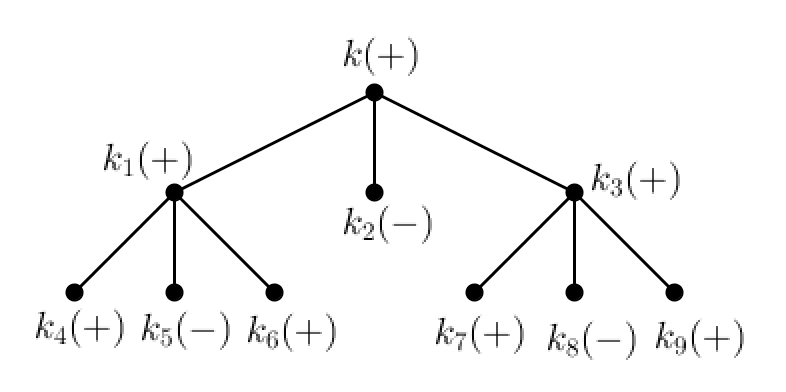}
  \caption{A tree $\Tc$ of order $3$ with its signature and decoration by wave numbers $k_j$.}\label{fig:tree}
\end{figure}

Diagrammatic expansions in terms of Feynman diagrams are the standard tool for the derivation of kinetic models. These are merely representations of the Picard (or Duhamel) iterates of the solution to the microscopic system, which is (\ref{NLS}) in this case. More precisely, we expand the solution to \eqref{NLS} in Fourier space as 
\begin{equation}\label{Picexp}
A_k(t)=A_k^{(0)}(t)+A_k^{(1)}(t)+\ldots+A_{k}^{(N)}+\textrm{remainder}
\end{equation}
where $A_k^{(n)}$ is the $n$-th Picard iterate of \eqref{Feqn}. This iterate itself is a sum over ternary trees with $n$ branching nodes, each depicting an interaction history affecting the mode $A_k$. We say that $n$ is the order of the tree. An example of such diagram is given in Figure \ref{fig:tree} in which the tree $\Tc$ is decorated at its nodes by wave numbers $k_j \in \Zb^d_L$. The root is always decorated by $k$ and for each branching node $\nf$ with children $\nf_1, \nf_2, \nf_3$ from left to right, there holds that $k_{\nf}=k_{\nf_1}-k_{\nf_2}+k_{\nf_3}$ where $k_{\nf}, k_{\nf_1},k_{\nf_2}, k_{\nf_3}$ are the decorations of $\nf,\nf_1, \nf_2, \nf_3$. We also associate a natural signature to the tree as follows: The signature of the tree is the same as that of its root node, and for each branching node $\nf$ with signature $\zeta_\nf$, the signature of its left and right child are $\zeta_\nf$ whereas that of the middle child is $-\zeta_\nf$. 

Up to phase shifting and rescaling to the kinetic times $t= \tau \cdot \delta T_{\mathrm{kin}}=\tau \cdot \delta L^{2\gamma}$ with $0\leq \tau \leq 1$, the contribution of this tree $\Tc$ to $A_k^{(n)}$, denoted by $(\Jc_{\Tc})_k(\tau)$, is given by
\begin{equation}\label{defjt}(\Jc_\Tc)_k(\tau)=\bigg(\frac{\delta}{2L^{d-\gamma}}\bigg)^n\widetilde{\zeta}(\Tc)\sum_{\Ds}\int_\Dc\prod_{\nf\in\Nc}e^{\pi i\zeta_\nf\delta L^{2\gamma}\Omega_\nf t_\nf}\,\mathrm{d}t_\nf\prod_{\lf\in\Lc}\sqrt{n_{\mathrm{in}}(k_\lf)}\cdot\eta_{k_\lf}^{\zeta_\lf}(\omega).
\end{equation} Here in (\ref{defjt}), the sum is taken over the set $\Ds$ of all decorations (or choices of wave numbers), $\Nc$ is the set of branching nodes of the tree $\Tc$, $\Lc$ is the set of leaves, $\widetilde{\zeta}(\Tc)=\prod_{\nf\in\Nc}(i\zeta_\nf)$, and $\Dc$ is the domain of time integration
\begin{equation}\label{defdomain}\Dc=\left\{(t_\nf)_{\nf \in \Nc}:0<t_{\nf'}<t_\nf<\tau,\mathrm{\ whenever\ \nf'\ is\ a\ child\ of\ \nf}\right\}.
\end{equation}

\begin{figure}[h!]
  \includegraphics[scale=0.45]{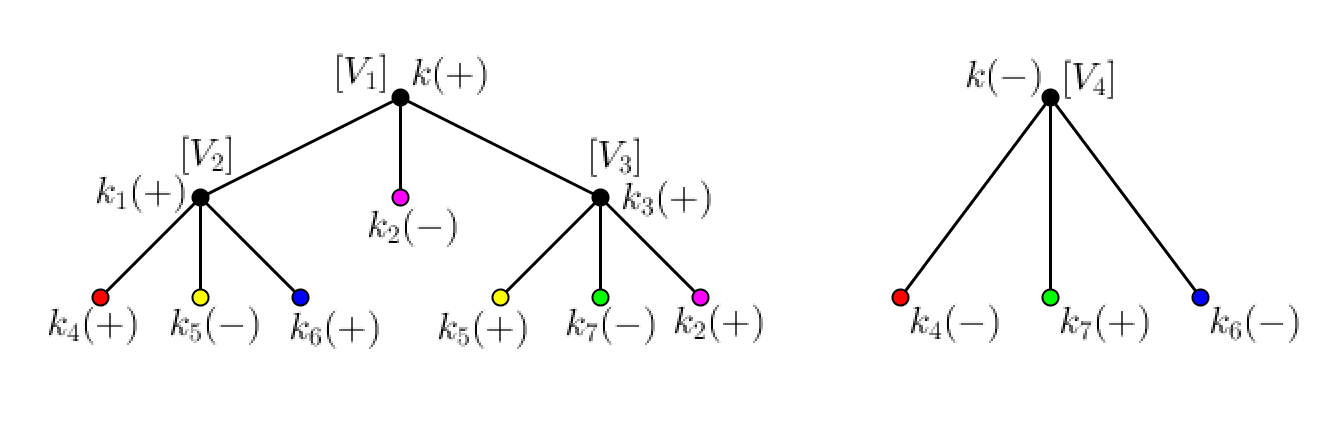}
  \caption{A couple of paired trees, shortly called a \emph{couple}, where each leaf is paired to exactly one other leaf (possibly from the same tree). Paired leaves are decorated with the same wave number. Roots have same decoration but opposite signatures.}
  \label{fig:Cpl}
\end{figure}

As such, when we study the second moments $\Eb |A_k|^2$ we end up with the correlation $\Eb (\Jc_{\Tc_1}\overline {\Jc_{\Tc_2}})$ of such expressions corresponding to trees $\Tc_1$ and $\Tc_2$. An inspection of \eqref{defjt} shows that this correlation involves taking the correlations $\Eb( \prod_{\lf\in\Lc_1 \cup \Lc_2} \eta_{k_\lf}^{\zeta_\lf})$ where $\Lc_1$ and $\Lc_2$ are the leaf-sets of the trees $\Tc_1$ and $\Tc_2$. Due to independence and mean zero assumptions on $\eta_k$, this expression is zero unless there is a complete pairing between the leaves in $\Lc_1\cup \Lc_2$ such that each positive leaf is paired to a negative one, and each such pair of leaves is decorated by the same wave number. This introduces the fundamental structure in our analysis, namely \emph{couples}, which are pairs of trees $\Tc^+$ and $\Tc^-$ of opposite signs such that each leaf in the couple is paired to exactly one leaf of opposite sign\footnote{These are the only structures that emerge in the study of second moments when $\eta_k$ are Gaussians, as a consequence of the Wick-Isserlis Theorem. In the non-Gaussian case, considered in \cite{DH21-2}, one has to also take account of couples where leaves are paired in sets of four or more. But these are shown to have lower order contributions to the second moments. The difference between Gaussian and non-Gaussian statistics starts to appear in moments $2p \geq 4$, for which one has to introduce a generalization of couples, namely gardens, which are $2p$-tuples of trees with their leaves paired to each other. We refer to \cite{DH21-2} for details on this.}. See Figure \ref{fig:Cpl}.

\subsection{Leading diagrams}
\begin{figure}[h!]
  \includegraphics[scale=0.16]{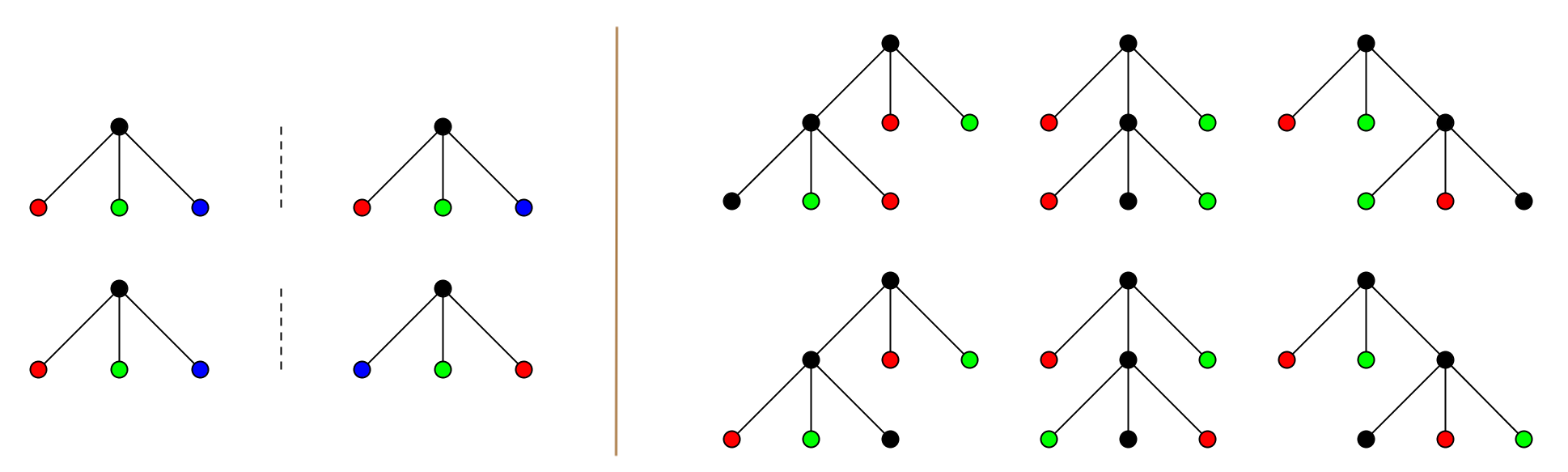}
  \caption{On the left are the two $(1,1)$ mini-couples, and on the right are the trees of order $2$ (two branching nodes) that are paired with the trivial tree to form the $(2,0)$ mini-couples.}
  \label{fig:minicouples}
\end{figure}

Recall the quantity computed in \eqref{formal1}, which is nothing but the correlation $\Eb[ (\Jc_{\Tc_1})_k \overline{(\Jc_{\Tc_1})_k}]$ where $\Tc_1$ is the (only) ternary tree with one branching node. This product corresponds to the so-called \emph{$(1,1)$ mini-couple} depicted on the left of Figure \ref{fig:minicouples}, and its contribution, in the kinetic limit, is the first Picard iterate of one of the four terms appearing in \eqref{wkenon}. The first iterate of the other three terms come from the contribution of the so-called $(2,0)$ mini-couples in which the trees on the right of Figure \ref{fig:minicouples} are paired with the trivial tree with one single node. As such, the first iterate of \eqref{wke} can be recovered as the kinetic limit of the contribution of the ($(1,1)$ and $(2, 0)$) mini-couples shown in Figure \ref{fig:minicouples}.

\FloatBarrier

\begin{figure}[h!]
  \includegraphics[scale=0.23]{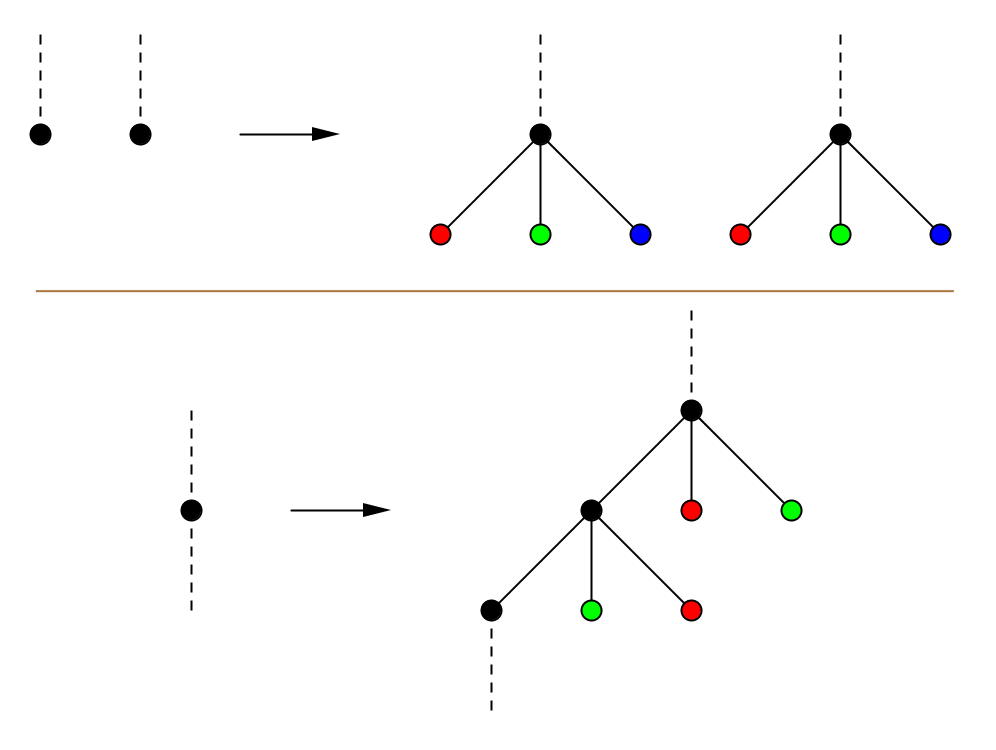}
  \caption{Regular couples are formed by starting with a trivial couple of only two paired leaves, and inductively applying the two operations depicted in the picture. The first operation (on top) replaces a pair of leaves by a $(1,1)$ mini-couple. The second operation (on bottom) replaces any node in the tree by one of the order-2 trees depicted on the right of Figure \ref{fig:minicouples}.}
  \label{fig:RegForm}
\end{figure}

But if this is the case, then the higher Picard iterates of \eqref{wke} have to come from the contribution of couples formed solely by concatenating the above mini-couples in the natural way, i.e. couples formed by mini-couples as building blocks. This generates what we call \emph{regular couples} which can be defined inductively as shown in Figure \ref{fig:RegForm}.

The regular couples play the role of the leading diagrams in the expansion, in the sense that they are also the ones which naturally optimize the counting estimates on couples. Indeed, by looking at \eqref{defjt}, one notices that the time-dependent part of the sum is given solely by the integral over the set $\Dc$ of the oscillating factors $e^{\pi i\zeta_\nf\delta L^{2\gamma}\Omega_\nf t_\nf}$. A naive bound on this time integral gives a function that barely fails to be in $L^1$ in $\Omega_{\nf}$. If we assume for the sake of argument that this function is indeed $L^1$, then one can use this summability to fix the values of $\lfloor\delta L^{2\gamma}\Omega_{\nf}\rfloor$, and reduce the estimate of $\Eb (\Jc_{\Tc^+}\overline {\Jc_{\Tc^-}})$ to a counting problem for the couple $\Qc=(\Tc^+, \Tc^-)$, where one counts the number of decorations of the couple subject to the extra conditions that $\delta L^{2\gamma}\Omega_{\nf}$ belongs to a fixed unit interval. 

As such, another way to define regular couples is the following: they are exactly the couples for which the corresponding counting problem decouples into ``effectively independent" three-vector counting problems of form
$$
a-b+c=h, \qquad |a|^2-|b|^2+|c|^2-|h|^2\in [M,M+\delta^{-1} L^{-2\gamma}],
$$ 
where $h\in \Zb^d_L$ and $M\in \Rb$ are fixed parameters and $(a,b,c)$ are variables in the counting problem. Indeed, for a regular couple with order $n$ (which is defined as the total number of branching nodes in both trees), one gets $n/2$ such counting problems, each amounting to the factor $\delta^{-1}L^{2d-2\gamma}$. Given the prefactor $\delta^n L^{-(d-\gamma)n}$ in \eqref{defjt}, one obtains that the contribution of such regular couples is $\delta^{n/2}$.

\subsection{Criticality and factorial divergence}

The goal is now to show that the contribution of all the remaining couples is an error term compared to that of regular couples. However, in this task we are faced with two major challenges. The first one is the \emph{probabilistic criticality} of the problem (see \cite{DNY3}), and the second is \emph{factorial divergence}. The criticality can be seen from the above estimate, namely that regular couples of order $n$ have contribution $O(\delta^{n/2})$. This estimate only improves by factors of $\sqrt \delta=O(1)$ with the order $n$ of the couple, rather than by a quantity that is $o(1)$ in the kinetic limit. Such diagrammatic expansions are notoriously difficult to treat mathematically, and as far as we know, the above theorem seems to be the first to treat a nonlinear out-of-equilibrium problem with such probabilistically critical nature. 

The second major difficulty of factorial divergence comes from the fact that there are factorially many couples of order $n$. If the only bound one can prove on such couples is a uniform upper bound like $\delta^{n/2}$ (which is sharp for regular couples), then one would be left with a series whose terms are bounded by $\delta^{n/2}n!$, which is a uselessly divergent bound. Indeed, overcoming this factorial divergence will be the main objective of the rigidity strategy discussed in Section \ref{section:rigidity} below.

\subsection{Molecules} The estimate above for regular couples is reduced to a counting problem for couples, namely counting the number of possible decorations by wave numbers satisfying many interdependent conditions. This counting problem turns out to be central to the estimation of the remaining couples as well. While the counting can be done more or less easily for regular couples, it gets quite more complicated for more general ones, for which one has to devise effective schemes and algorithms. One of the main tools introduced in \cite{DH21} to perform this counting procedure effectively is the notion of ``molecules". Each couple, like the one depicted in Figure \ref{fig:Cpl}, is associated to exactly one molecule which is a directed graph that carries all the needed combinatorial information from the couple, but in a much more effective way. The molecule associated to the couple in Figure \ref{fig:Cpl} is shown in Figure \ref{fig:Mole} below.

\begin{figure}[h!]
  \includegraphics[scale=0.5]{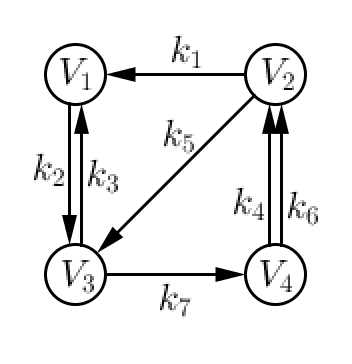}
  \caption{The molecule corresponding to the couple in Figure \ref{fig:Cpl}. Each atom in the molecule corresponds to a branching node in the couple, as labeled. Two such atoms are linked if their corresponding nodes in the couple either a) are children of each other, or b) have children that are leaves paired to each other.}
  \label{fig:Mole}
\end{figure} 

Each vertex of the molecule, called an \emph{atom}, corresponds to a branching node in the couple (see for example Figure \ref{fig:Cpl}), and each edge, called a \emph{bond}, between the atoms signifies either a parent-child relationship between the two branching nodes of the couple or a leaf pairing between their children. The directions of the bonds carry the signature relations in the couples. With this new combinatorial object, the counting problem can now to reformulated as an algorithm that ``digests" this molecule by breaking its bonds until we are left with several much more tractable and independent counting problems. This the essence of the combinatorial arguments needed to execute the strategy discussed below (cf. Section \ref{section:rigidity}).

\subsection{Diagrammatic Cancellations and scaling laws} \label{cancellations}
\begin{figure}[h!]
 \hspace*{-0.5cm} \includegraphics[scale=0.2]{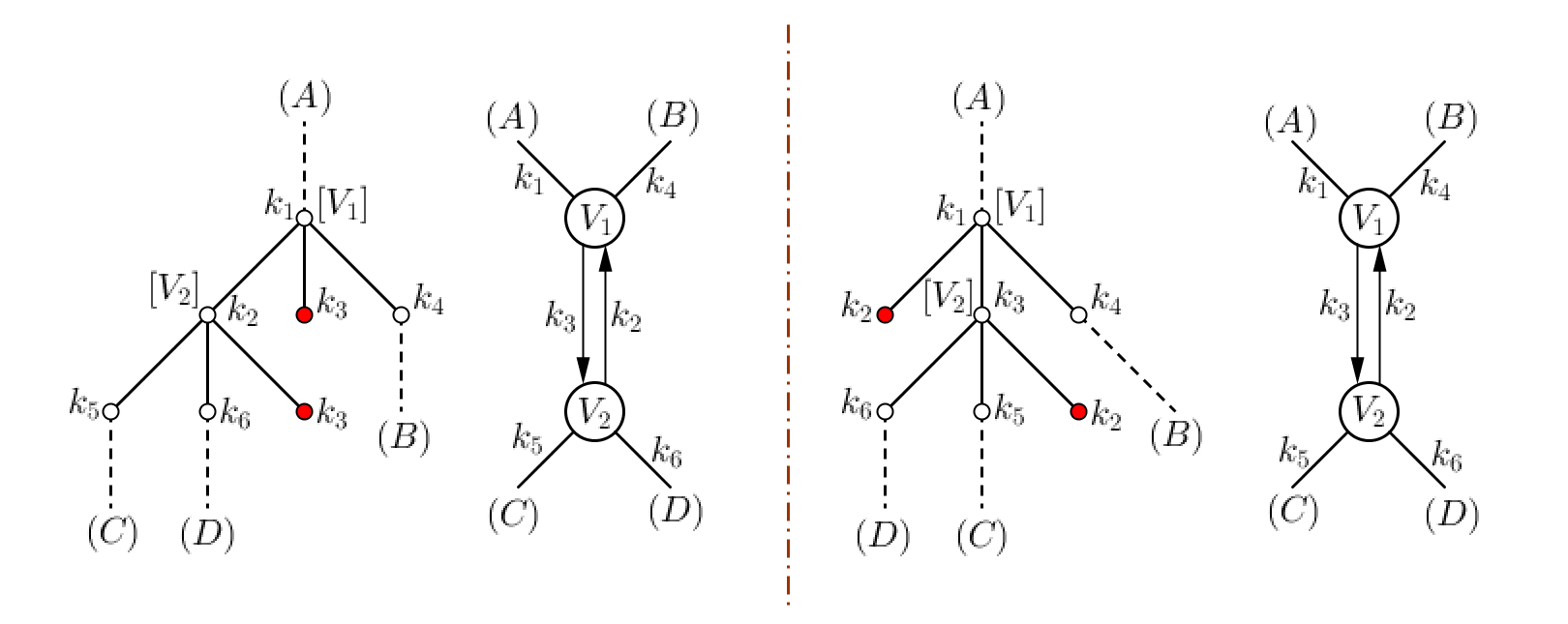}
  \caption{A cancellation happens between the two couples depicted on the left and the right. The couples are drawn next to their respective molecules. One notices that those two cancelling couples have the same molecule. }
  \label{fig:Cancel1}
\end{figure}

Unfortunately, the above two difficulties (criticality and factorial divergence) are not the only ones to proving the convergence of the diagrammatic expansion. As was first uncovered in our first joint work \cite{DH19} on the subject, there exists a family of non-regular couples whose contribution is almost just as large as that of regular couple. We called such couples \emph{irregular chains} in \cite{DH21} where we treated scaling laws near $\gamma=1$. The only possible way to prove that such diagrams do not disturb the leading regular ones in the kinetic limit, is to show that they cancel in some form. This is one of the main novelties of \cite{DH21} where a major cancellation is exhibited at the top order for such diagrams, in a way that shows their contribution to be sub-leading in the end. See Figure \ref{fig:Cancel1}.

\begin{figure}[h!]
\hspace*{-0.5cm} \includegraphics[scale=0.4]{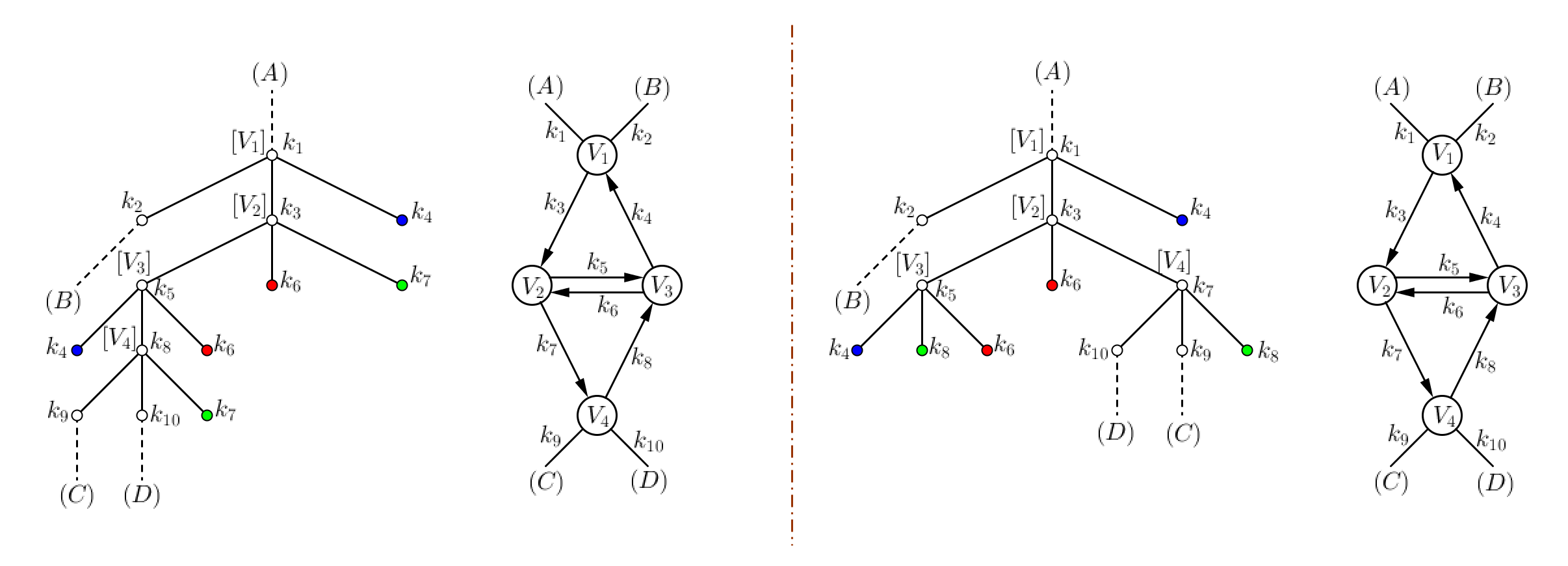}
  \caption{The much more elaborate cancellation that is needed for scaling laws $\gamma<\frac23$. Again in this case, the couples that cancel have the same molecule.}
  \label{fig:Cancel2}
\end{figure}

Even more surprisingly, other, much more complicated, large diagrams emerge when treating the other scaling laws $\gamma\in (0,1)$, which is done in our most recent work \cite{DH22}. These diagrams are lower order error terms for the scaling laws with $\gamma\approx 1$, so they did not appear in our work \cite{DH21}, but they do become large for scaling laws with $\gamma< \frac 23$. These diagrams get most elaborate for the scaling law $\gamma=\frac12$, so much so that it becomes hard to identify them at the level of couples. There, the right place to identify such large structures is by looking at the so-called ``molecule picture" associated to each couple. Indeed, the algorithm that resolves the counting problem associated to the molecules identifies the large ``divergent" ones for which there is no gain compared to the leading regular couples. In a way, such divergent molecules, and the couples they represent, cannot be summed absolutely in the expansion. We show in Figure \ref{fig:Cancel2} pictures of such divergent couples alongside the molecules corresponding to them. 

The miraculous fact is that these structures also come in pairs, just like the irregular chains mentioned above, although the exact ``twist" going from one to the other becomes less trivial to identify. These pairs cancel at top order, which allows us to show that their contribution does not disturb the leading contribution of regular couples in the kinetic limit. The exact physical meaning of such diagrams and their cancellation is still not clear; in fact, the large size and complexity of this cancellation seem to be novel even in the physics  literature. 

\subsection{Rigidity Theorem}\label{section:rigidity}

At this point, we have discussed the largest diagrams appearing in the expansion: some of those are the leading regular couples which give the needed kinetic limit, whereas others are highly complex diagrams that are large but exhibit elaborate cancellations between themselves\footnote{Luckily, the number of regular diagrams as well as those needing cancellation grows exponentially rather than factorially in the order $n$ of the couple.}. What about the rest? Recall here that the main issue is the factorial divergence in the number of couples of order $n$. Because of this, it is not enough just to show that each non-leading diagram of order $n$ gains a constant power compared to the leading diagrams (regular couples) of order $n$; In fact, these power gains should be strong enough to offset the factorial number of these diagrams. 

The strategy is to classify couples of order $n$ according to some index $1\leq r\leq n$, such that a) the number of couples of index $r$ is $\lesssim C^n r!$, and b) the estimate for couples of index $r$ gains a power comparable to $L^{-cr}$. This gain is enough to offset the factorial divergence since $C^n r! L^{-cr}=o(1)$ provided we only expand in \eqref{Picexp} up to $N\sim \log L$ (recall that $r\leq n \leq N$). This index $r$ turns out to be basically the order of the reduced couple obtained from the original one after performing sort of a surgery that excises out the regular sub-couples in the couple, as well as all the other large structures mentioned in Section \ref{cancellations}. The fact that all these structures have exponential (rather than factorial) complexity allows to bound the number of couples with the same reduced couple by $C^n$, which gives the estimate $C^n r!$ for the number of couples with index $r$. The execution of this strategy is referred to as a ``Rigidity Theorem" in \cite{DH21, DH22}, which is probably the most important theorem in those works.

Finally, the estimate on the remainder in the expansion \eqref{Picexp} is based on obtaining good resolvent bounds for the linearized operator that appears in the equation for the remainder. Such estimates turn out to follow from similar arguments to those involved in the analysis of the diagrams themselves. 
\section{Conclusion}
We have laid down the mathematical theory of wave turbulence for NLS for the range of scaling laws $\gamma\in (0,1]$ of interest to both mathematicians and physicists. This range is the most physically relevant one as it does not require any conditions on the shape of the box (for $\gamma<1$). The importance of introducing and clarifying the notion of scaling laws, both from mathematical and physical perspectives, have been discussed. A delicate diagrammatic analysis is needed for the mathematical theory, and one of its interesting consequences is the dependence on the scaling laws in identifying the largest diagrams that require subtle cancellations. These cancellation structures are still quite obscure, and their emergence seems like a mathematical miracle, so a better physical understanding of them is needed. 

\end{document}